# Plasma Shock Layer Equations

Y. Farjami

Departmentof Computer & IT,

University of Qom,

Qom, Iran.

Email: farjami@qom.ac.ir

## 1  Introduction

Many of physical environments consist of two or more types of particles, each having its own characteristics, such an environment may undergo some evolutions in the presence of external forces, e.g., in the presence of electromagnetic or gravitation fields, and may also be affected by a source of energy or by interaction with its external environment. The situation becomes more complex when there are some internal interaction due to collisions between particles, chemical reactions and/or various dissipation mechanisms. From microscopic point of view, the random motion of different particles and local forces, energy or mass transaction and their recursive dependence with that of global makes the qualitative behaviour of the environment even more complex. Heat convection in a liquid mixture, a vibrating or rotating viscose fluid, chemically reacting multiphase media and plasma in a electric and magnetic field are examples, to say a few.

Up to the authors' knowledge there is no general and precise mathematical formulation ef-ficiently describing such environments, although statistical formulation in Kinetic theory combined with some physical interpretations is a good approach but insufficient, say in describing radiation.

Our aim in this note is to give a set of partial differential equation governing evaluation of a multispecies environment consisting of charged particles in an electromagnetic field.

## 2 Fundamental Principles and Equations

Conservation laws for a physical system are usually derived from some physical principles such as conservation of mass, momentum, energy and charge. We are dealing with a system consisting of $n$ species, such a system is called a plasma of $n$ component. Each species or component will be identified by a subscript $\alpha$.

We use the following notations for physical or thermodynamical quantities.

$m$ mass
$u = (u^1, u^2, u^3)$ velocity vector
$x = (x^1, x^2, x^3)$ coordinate of a point
$t$     time
$\tau^{ij}$     stress tensor
$p$     pressure
$E$     Electric field
$B$     magnetic field density
$J$     current density
$T$     Temperature
$Q$     Heat flux
$p$     density
$\rho_e$     charge density
$n$     number of particle in unit volume.

Thus with these notations $m_\alpha$ denotes the mass of a particle of $\alpha$-species in a unit volume. We also use the convention that the $i$-th component of a vector say $u_\alpha$ is denoted by $u_\alpha^i$, and that the repeated superscripts denote summation, e.g. $u_\alpha^i u_\beta^i = \sum_{i=1}^{3} u_\alpha^i u_\beta^i = u_\alpha \cdot u_\beta$.

## 2.1 Equation of State.

For each species of plasma, there is a functional relation between $p_\alpha, \rho_\alpha$ and $T_\alpha$, i.e.,

$$p_\alpha = p_\alpha(\rho_\alpha, T_\alpha). \tag{1}$$

## 2.2 Equation of Continuity.

The conservation of mass of $\alpha$-th species given the equation of continuity of this species which may be written as

$$\frac{\partial \rho_\alpha}{\partial t} + \frac{\partial}{\partial x^i}(\rho_\alpha u^i_\alpha) = \sigma_\alpha. \tag{2}$$

Here $\sigma_\alpha$ is the mass source per unit volume of $\alpha$-th species which may be due to ionization processes or other chemical reactions. If there is no chemical reaction and the degree of ionization is constant, $\sigma_\alpha$ is zero. Even if some of the $\sigma_\alpha$ are not equal to zero, by conservation of mass of the plasma as a whole we have

$$\sum_{\alpha=1}^{n} \sigma_\alpha = 0. \tag{3}$$

## 2.3 Equation of Motion.

Considering the conservation of momentum for $\alpha-$th species, then we get the equation of motion for $\alpha-$th species as follows:

$$\frac{\partial \rho_\alpha u^i_\alpha}{\partial t} + \frac{\partial}{\partial x^j}(\rho_\alpha u^i_\alpha u^j_\alpha - \tau^{ij}_\alpha) = X^i_\alpha + \sigma_\alpha + \sigma_\alpha Z^i_\alpha \tag{4}$$

where $\sigma_\alpha Z^i_\alpha$ is the $i-$th component of the momentum source per unit volume associated with the mass source $\sigma_\alpha$. We also have

$$\sum_{\alpha-1}^{n} \sigma_\alpha Z^i_\alpha = 0. \tag{5}$$

The term $\tau^{ij}_\alpha$ is the $ij-$th component of the stress tensor of $\alpha-$th species. In macroscopic analysis we can only postulate the relation between $\tau^{ij}_\alpha$ and other physical quantities. For instance, the simple relation is that $\tau^{ij}_\alpha$ consists of total pressure $p_{t\alpha}$ and viscouse stress tensor $\tau^{ij}_{\alpha 0}$ which may be assumed to be proportional to velocity gradient according to Naiver-Stokes relation. Of course the stress $\tau^{ij}_\alpha$ may also depend on temperature gradient or on concentration gradient. The total pressure $p_{ts}$ includes both gasdynamic pressure $p_\alpha$ and the radiation pressure $p_{R\alpha}$. At low or moderately high temperature, the radiation pressure is usually negligible. Only at very high temperature the radiation pressure should be included in the analysis of flow problem.

The term $X_\alpha^i$ is the $i-$th component of the body force per unit volume. It consists of the electromagnetic force and nonelectric body force. We may write

$$X_\alpha^i = F_{e\alpha}^i + F_{g\alpha}^i \tag{6}$$

where $F_{g\alpha}^i$ is the $i-$th component of nonelectric forces such as gravitational force etc.

The electromagnetic force $F_{e\alpha}^i$ may be written in the following form

$$F_{e\alpha}^i = \rho_{e\alpha}[E^i + \mu_e(u_\alpha \times H)^i] + F_{e\alpha o}^i \tag{7}$$

here $\rho_{e\alpha}$ is the electric charge density of $\alpha$-species, hence $\rho_{e\alpha} = \nu_\alpha e_\alpha$ where $e_\alpha$ is the charge on a particle of $\alpha$-species and $\nu_\alpha$ is the number of particles of $\alpha$s per unit volume. The velocity vector of $\alpha$-species is $u_\alpha$. Finally $F_{e\alpha 0}$ is the interaction force on the $\alpha$-species due to all other kind of particles in the plasma. By Newton's third law of motion, we have

$$\sum_{\alpha=1}^{n} F_{e\alpha o}^i = 0. \tag{8}$$

The term $H$ is the vector of magnetic field strength.

## 2.4 Equation of Energy.

The conservation of energy of $\alpha$-species gives the equation of energy of this species which may be written as

$$\frac{\partial \bar{e}_\alpha}{\partial t} + \frac{\partial}{\partial x^j}(\bar{e}_\alpha u_\alpha^j - u_\alpha^i \tau_\alpha^{ij} - Q_\alpha^j) = \varepsilon_\alpha \tag{9}$$

where $\bar{e}_\alpha$ = total energy of $\alpha$-species in the plasma per unit volume. we can write $\bar{e}_\alpha = \rho_\alpha \bar{e}_{m\alpha}$ where $\bar{e}_{m\alpha}$ is total energy of $\alpha$-species of the plasma per unit mass and is give by

$$\bar{e}_{m\alpha} = U_{m\alpha} + \frac{1}{2}u_\alpha^2 + \phi_\alpha + \frac{E_{R\alpha}}{\rho_\alpha}, \tag{10}$$

where

$U_{m\alpha}$ = internal energy of $\alpha$-species per unit mass,

$\frac{1}{2}u_\alpha^2$ = kinetic energy per unit mass of $\alpha$-species.

$\phi_\alpha$ = potential energy per unit mass of $\alpha$-species.

$E_{Rs}$ = radiation energy of $\alpha$-species per unit volume.

The first term on the left hand side of (2.9) is the rate of change of energy of $\alpha$-species per unit volume. The second term is the energy flow by convention, i.e., $(\frac{\partial}{\partial x^j}\bar{e}_\alpha u_\alpha^j)$. The third term $-(\frac{\partial}{\partial x^j} u_\alpha^j \tau_\alpha^{ij})$ is the rate of energy dissipation by stress tensor $\tau_\alpha^{ij}$. The fourth term $-(\frac{\partial}{\partial x^j}, Q_\alpha^j)$ is the energy change due to heat flow, which consists of both head conduction, $Q_{c\alpha}^j$ and radiation, $Q_{R\alpha}^j$; i.e.,

$$Q_\alpha^j = Q_{c\alpha}^j + Q_{R\alpha}^j. \tag{11}$$

The term on the rhs of equation (2.9) is the energy production term or energy source per unit volume. It consists of terms due to electromagnetic field $\varepsilon_{e\alpha}$ and due to chemical reaction $\varepsilon_{c\alpha}$.

## 2.5 Maxwell's Equation of electromagnetic Field.

The equations governing the electric and magnetic fields are the Maxwell's equations which are

$$\nabla \times H = J + \frac{\partial D}{\partial t}, \tag{12}$$

$$\nabla \times E = -\frac{\partial B}{\partial t}, \tag{13}$$

$$D = \varepsilon E, \tag{14}$$

$$B = \mu_e H. \tag{15}$$

Here $H$ is the magnetic field strength, $E$ is the electric field strength, $J$ is the electric current density, $\varepsilon$ is the inductive capacity or dielectric constant, $B$ is the magnetic flux density and $\mu_e$ is the magnetic permeability. Ordinarily we may assume that both $\varepsilon$ and $\mu_e$ are constant.

The interaction of electromagnetic field equations with the flow of plasma is through the term $J$ which is

$$J^i = \sum_{\alpha=1}^{n} \rho_{e\alpha} u_\alpha^i. \tag{16}$$

The electric current density $J$ is due to the motion of various charged particles in the plasma. Here we don't give relations for $J$ in terms of electromagnetic variables, such a relation which is usually referred to generalized ohm law will be studied is a subsequent section.

## 2.6 Mean variables of the plasma as a whole.

Even though it is sufficient to use the partial quantities $(\rho_\alpha, \rho_\alpha, T_\alpha, u_\alpha^i)$ to describe the motion of plasma, it is sometimes more convenient to use the mean quantities of the plasma as a whole in our analysis at plasma dynamics. The relations between the mean quantities and the corresponding partial quantities are defined as follows.

(a) Density $\rho$ of plasma.

The number density $\nu$ of the plasma is the sum of the number density $\nu_\alpha$ of all the species, i.e.,

$$\nu = \sum_{\alpha=1}^{n} \nu_\alpha. \tag{17}$$

The mass density or simply the density $\rho$ of the plasma is

$$\rho = m\nu = \sum_{\alpha=1}^{n} m_\alpha \nu_\alpha = \sum_{\alpha=1}^{n} \rho_\alpha, \tag{18}$$

where $m$ is the mean mass of a particle in the plasma.

(b) Pressure $p$ of the plasma.

The pressure $p$ of the plasma is the sum of all the partial pressures of the species $p_\alpha$ in the plasma, i.e.,

$$p = \sum_{\alpha=1}^{n} p_\alpha. \tag{19}$$

(c) Temperature $T$ of the plasma.

The kinetic temperature $T$ of the plasma is defined as follows

$$T = \frac{1}{\nu} \sum_{\alpha=1}^{n} \nu_\alpha T_\alpha \tag{20}$$

(d) Flow velocity $u^i$ of the plasma.

The $i$-th component of the flow velocity $u^i$ of the plasma as a whole is defined by the following relation

$$\rho u^i = \sum_{\alpha=1}^{n} \rho_\alpha u_\alpha^i. \tag{21}$$

The diffusion or random velocity $w_\alpha^i$ of the $\alpha$-species in the plasma is then

$$w_\alpha^i = u_\alpha^i - u^i. \tag{22}$$

It is evident that

$$\sum_{\alpha=1}^{n} \rho_\alpha w_\alpha^i = 0 \tag{23}$$

(e) Excess electric charge $\rho_e$.

The sums of the charge densities $\rho_{e\alpha}$ of all species gives the excess electric charge $\rho_e$ of the plasma, i.e.,

$$\rho_e = \sum_{\alpha=1}^{n} \rho_{e\alpha} = \sum_{\alpha=1}^{n} e_\alpha \nu_\alpha. \qquad (24)$$

(f) Electrical current density $J^i$.

The electrical current density in the plasma is defined as

$$J^i = \sum_{\alpha=1}^{n} \rho_{e\alpha} u_\alpha^i \sum_{\alpha=1}^{n} \rho_{e\alpha} w_\alpha^i + u^i \sum_{\alpha=1}^{n} \rho_{e\alpha} = \iota^i + \rho_e u^i, \qquad (25)$$

where $f\iota^i$ is the $i$th component of the conduction current and $\rho_e u^i$ is the $i$th component of the convection current.

## 2.7

Fundamental equations for the mean quantities of the plasma the fundamental equations for the mean quantities of the plasma may be derived from the equations for partial quantities. These equations for the mean quantities are given below.

(a) Equation of state of the plasma.

Similar to (2.1) we may assume a general form for equation of state of plasma, as follows

$$p = p(\rho, T), \qquad (26)$$

where $\rho$, $p$, and $T$ are defined by (2.18), (2.19) and (2.20) respectively.

(b) Equation of continuity of the plasma.

If we add all the $n$ equations of the type (2.2) and we the relation (2.3), we have

$$\frac{\partial \rho}{\partial t} \frac{\partial}{\partial x^i}(\rho u^i) = 0, \qquad (27)$$

which is called the equation continuity.

(c) The equation of motion of the plasma.

If we add all the $n$ equations of the type (2.4) and using (2.5), $w$ have

$$\frac{\partial \rho u^i}{\partial t} \frac{\partial}{\partial x^j}(\rho u^i u^j - \tau^{ij}) = X^i = F_e^i + F_g^i. \qquad (28)$$

This is the equation of motion of the plasma in the direction of $i$th axis. The various terms in equation (2.28) are explained as follows.

The term $F_g^i$ is the nonelectric force, i.e.,

$$F_g^i = \sum_{\alpha=1}^n F_{g\alpha}^i. \qquad (29)$$

The term $F_e^i$ is the electromagnetic force, i.e.,

$$F_e^i = \sum_{\alpha=1}^n F_{e\alpha}^i = \sum_{\alpha=1}^n \rho_{e\alpha}(E^i + \mu_e(u_\alpha \times H)^i) = \sum_{\alpha=1}^n \rho_{e\alpha}E^i + \mu_e(J_\alpha \times H)^i$$
$$= \rho_e E^i + \mu_e(J \times H)^i, \qquad (30)$$

where $J_\alpha = \rho_{e\alpha} u_\alpha$ is the electric current density of $\alpha$-species and hence by (2.25) we have

$$J = \sum_{\alpha=1}^n J_\alpha. \qquad (31)$$

The stress tensor $\tau_\alpha^{ij}$ of the $\alpha$-species may be written as

$$\tau_\alpha^{ij} = -\delta^{ij} p_{t\alpha} + \tau_{\alpha 0}^{ij} \qquad (32)$$

where $\delta^{ij}$ is the Kronecker delta and $\tau_{\alpha 0}^{ij}$ is the viscous stress tensor and $p_{t\alpha}$ is the total pressure of as defined by

$$p_{t\alpha} = p_\alpha + p_{R\alpha}, \qquad (33)$$

where $p_\alpha$ is the gas dynamic pressure and $p_{R\alpha}$ is the radiation pressure of $\alpha$-species. The total stress tensor $\tau^{ij}$ of the plasma is

$$\tau^{ij} = -\sum_{\alpha=1}^n \delta^{ij} p_{t\alpha} - \sum_{\alpha=1}^n \rho_\alpha w_\alpha^i w_\alpha^j + \sum_{\alpha=1}^n \tau_{\alpha 0}^{ij}$$
$$= -\delta^{ij} p_t + \tau^{ij}, \qquad (34)$$

where the total viscous stress tensor is

$$\tau_0^{ij} = \sum_{\alpha=1}^n \tau_{\alpha 0}^{ij} - \sum \rho_\alpha w_\alpha^i w_\alpha^j. \qquad (35)$$

The total viscous stress tensor is due to a very complicated molecular motion. For first approximation, we may assume that $\tau_0^{ij}$ is of the following form

$$\tau_0^{ij} = \mu\left(\frac{\partial u^i}{\partial x^j} + \frac{\partial u^j}{\partial x^i}\right) + \mu_1 \frac{\partial u^k}{\partial x^k} \delta^{ij}, \qquad (36)$$

where $\mu$ is the ordinary or first coefficient of viscosity and $\mu_1$ is the second coefficient of viscosity. Ordinarily we may assume that $\frac{2}{3}\mu + \mu_1 \geq 0$. Sometimes $\mu$ is called the bulk or longitudinal viscosity and $\mu_1$ is called the transverse or shear viscosity. The equation (2.36) may be written as

$$\tau_0^{ij} = \mu\left(\frac{\partial u^i}{\partial x^j} + \frac{\partial u^j}{\partial x^i} - \frac{2}{3}\frac{\partial u^k}{\partial x^k}\right) + \left(\frac{2}{3}\mu + \mu_1\right)\frac{\partial u^k}{\partial x^k}\delta^{ij}$$
$$\mu\left(\frac{\partial u^i}{\partial x^j} + \frac{\partial u^j}{\partial x^i} - \frac{2}{3}\frac{\partial u^k}{\partial x^k}\right) + \zeta\frac{\partial u^k}{\partial x^k}\delta^{ij} \tag{37}$$

where $\zeta = \frac{2}{3}\mu + \mu_1$ is nonnegative. The coefficients of viscosity are functions of temperature and the composition of the plasma. They are also functions of electric and magnetic fields.

The total pressure consists of the gasdynamic pressure $p$ and the radiation pressure $p_R$, i.e.,

$$p_t = \sum_{\alpha=1}^{n} p_{t\alpha} = \sum_{\alpha=1}^{n} p_\alpha + \sum_{\alpha=1}^{n} p_{R\alpha}$$
$$= p + p_R. \tag{38}$$

The gasdynamic pressure is the value of pressure used in the equation of state. The radiation pressure is determined by a very complicated radiation phenomena. However for simplicity we may assume that

$$p_R = \frac{a_R}{3}T^4, \tag{39}$$

where $a_R$ is known as Stefan-Boltzmann constant.

(d) Equation of energy.

If we add all the $n$ equations of the type (2.9) we have

$$\frac{\partial \rho\overline{e_m}}{\partial t} + \frac{\partial}{\partial x^j}(\rho\overline{e_m}u^j - u^i\tau^{ij} - Q^j) = \varepsilon_t. \tag{40}$$

This is the energy equation of the plasma. The various terms in equation (2.40) are explained as follow.

The total energy of the plasma per unit mass is

$$\overline{e_m} = U_m + \frac{1}{2}u^2 + \phi + \frac{E_R}{\rho}, \tag{41}$$

where $U_m = \frac{1}{\rho}\sum_{\alpha=1}^{n}(\rho_\alpha U_{m\alpha} + \frac{1}{2}\rho_\alpha w_\alpha^2)$ is total energy of the plasma per unit mass.

$\phi = \frac{1}{\rho}\sum_{\alpha=1}^{n}\phi_\alpha\rho_\alpha$ is the total potential energy of the plasma per unit mass. $E_R = \sum_{\alpha=1}^{n} E_{R\alpha}$ is the total radiation energy of the plasma per unit volume. As a first approximation which is consistent with equation (2.39), we have

$$E_R = a_R T^4. \tag{42}$$

The heat flow $Q^j$ of the plasma consists of both heat conduction $Q_c^i$ and radiation $Q_R^j$, i.e.,

$$Q^j = Q_c^j + Q_R^j, \qquad (43)$$

where $Q_c^j$ the heat conduction flux, is

$$Q_c^j = \sum_{\alpha=1}^{n} \{Q_{c\alpha}^j + u^i \rho_\alpha w_\alpha^i w_\alpha^j - \rho_\alpha (U_{m\alpha} + \frac{1}{2} u_\alpha^2 + \phi_\alpha) w_\alpha^j + \tau_{\alpha 0}^{ij} w_\alpha^i - \rho_\alpha w_\alpha^j\}. \qquad (44)$$

The heat conduction flux is also due to a very complicated molecular motion. For first approximation, we may assume that

$$Q_c^j = k \frac{\partial T}{\partial x^j}, \qquad (45)$$

where $k$ is the coefficient of heat conductivity which is a function of temperature and composition of the plasma. It may also depend on the electric and magnetic fields.

The radiation flux is due to complex radiation phenomena. For a first approximation, we may assume that

$$Q_R^j = D_R \frac{\partial E_R}{\partial x^j}, \qquad (46)$$

where $D_R$ is known as the diffusion coefficient of radiation and may be written as

$$D_R = \frac{c l_R}{3} \qquad (47)$$

$$l_R = \frac{1}{K_R \rho} \qquad (48)$$

where $l_R$ is the Roseland mean free path of radiation and $K_R$ is the capacity or Roseland mean absorption coefficient. $K_R$ depends on both temperature $T$ and density $\rho$. Here $c$ is the velocity of light.

The energy source $\varepsilon_t$ consists of one due to chemical reaction $\varepsilon_c$ and one due to electromagnetic fields $\varepsilon_e$, i.e.,

$$\varepsilon_t = \sum_{\alpha=1}^{n} \varepsilon_\alpha = \sum_{\alpha=1}^{n} (\varepsilon_{\alpha c} + \varepsilon_{\alpha e}) = \varepsilon_c + \varepsilon_e, \qquad (49)$$

and

$$\varepsilon_e = \sum_{\alpha=1}^{n} \varepsilon_{\alpha e} = \sum_{\alpha=1}^{n} J_\alpha \cdot E = (\sum_{\alpha=1}^{n} J_\alpha) \cdot E = J^j \cdot E^j. \qquad (50)$$

(e) Equation of conservation of electrical charge.

If we multiply the factor

$$\gamma_\alpha = e_\alpha / m_\alpha, \qquad (51)$$

to equation (2.2), we have the equation of conservation of charge of each charge species which gives

$$\frac{\partial \rho_e}{\partial t} + \frac{\partial}{\partial x^j}(\rho_{ex} u_\alpha^j) = \gamma_\alpha \sigma_\alpha, \tag{52}$$

where we have used the filling simple relations.

$$\rho_{e\alpha} = \nu_\alpha e_\alpha = \nu_\alpha m_\alpha e_\alpha / m_\alpha = \gamma_\alpha \rho_\alpha. \tag{53}$$

The conservation of total electrical charge of the plasma gives

$$\sum_{\alpha=1}^{n} \gamma_\alpha \sigma_\alpha = 0. \tag{54}$$

If we add all the $n$ equations of type (2.53), we have the equation of conservation of electrical charge of the plasma or electrical equation of continuity of plasma, which is

$$\frac{\partial \rho_e}{\partial t} + \frac{\partial J^j}{\partial x^j} = 0. \tag{55}$$

(f) Equation electrical current.

Equation (2.4) multiplied by $\gamma_\alpha$ gives

$$\frac{\partial \rho_{e\alpha} u_\alpha^i}{\partial t} + \frac{\partial}{\partial x^j}(\rho_{e\alpha} u_\alpha^i u_\alpha^j - \gamma_\alpha \tau_\alpha^{ij}) = \gamma_\alpha (X_\alpha^i + \sigma_\alpha Z_\alpha^i). \tag{56}$$

If we add all the $n$ equations of type (2.26), we have

$$\frac{\partial J^i}{\partial t} + \frac{\partial}{\partial x^j} \sum_{\alpha=1}^{n}(\rho_{e\alpha} u_\alpha^i u_\alpha^j - \gamma_\alpha \tau_\alpha^{ij}) = \sum_{\alpha=1}^{n} \gamma_\alpha (X_\alpha^i + \sigma_\alpha Z_\alpha^i). \tag{57}$$

Equation (2.57) is the differential equation for the current density $J^i$. In general equation (2.57) is too complicated to be analyzed. Usually in the analysis of plasma dynamics we use generalized Ohm's law which will be derived in the next section. From macroscopic point of view when we consider the plasma as a whole the following from of Ohm's law will be sufficient.

$$i^i = \sigma[E^i + \mu_e(u \times H)^i], \tag{58}$$

where $\sigma$ is know as the electrical conductivity.

## 2.8 Fundamental equations of plasma dynamics.

There are several ways to describe the motion of plasma.

(a) Kinetic theoretic approach.

Since the plasma is composed of a large number of particles, charged and neutral, the most accurate but a more complicated method of description of the flow of a plasma is the molecular theory of plasma. However, because of many physical and mathematical difficulties, it is not possible at the present to treat the gas or plasma flow problem exactly by molecular theory. Many simplifications about the molecular forces and collisions have to be made in the formulation of theory and the resultant equations can only be solved approximately. Such a simplified molecular theory is known as the kinetic theory of gases. The generalization of the ordinary kinetic theory of gases to the case of plasma is one of the most interesting current research problems. Because of the importance of this subject to plasma dynamics, the reader should have a general idea of such a theory. Here we should mention that most of the previous equations can be derived using kinetic theory of plasma see [Akiezer, Wood]. For further information in kinetic theory the reader is referee to [Wood, Akheizer, Sutton].

(b) Multifluid description of plasma.

We may use the partial variables $(T_\alpha, p_\alpha, \rho_\alpha, u_\alpha^i, E^i, H^i)$ to describe the motion and the states of the plasma where $\alpha = 1, 2, \cdots, n$ and $i = 1, 2, 3$. The fundamental equations in this case are equation (2.1), (2.2), (2.4), (2.9), (2.10), and (2.11). There are $6n + 6$ variables and $6n + 6$ equations.

(c) Mean variable description.

We may use the mean variable $(T, p, \rho, u^i, E^i, H^i)$ and partial variable, $(T_\alpha, p_\alpha, \rho_\alpha, u_\alpha^i)$ where $\alpha = 1, 2, \cdots, n-1$. The fundamental equations are (2.10), (2.11), (2.19), (2.27), (2.28), (2.40), and the corresponding equations for the partial variables. We may also replace the variables $u_\alpha^i$ by $w_\alpha^i$ and replace (2.4) by a corresponding equation of $w_\alpha^i$. This approach is of great interest because the fundamental equations of magnetogasdynamics are derived from this approach.

(c) Fully ionized plasma.

As an special case of the proceeding approach let's consider a fully ionized plasma, in which there are only two type of particles, the electrons and ions, thus we have $n = 2$. We have the 18 variables $T, p, \rho, u^i, E^i, H^i, T_1, p_1, \rho_1, w_1^i$, where subscript 1 represents the values of one of the species in the plasma, e.g., those of the electrons.

(d) Single fluid description of plasma.

In many practical problems we may consider the plasma as a single fluid. This approach has been widely accepted and adapted by many authors. This description of plasma is usually referred to as MHD. In this case, our variables are quite similar to those of fully ionized plasma, except that we are not interested in the partial values $T_1$ and $p_1$. There are now only 16 variables in our problem., i.e., $T, p, \rho, u^i, E^i, H^i, \rho_e$ and $J^i$. The fundamental equations are then

$$\nabla \times H = J + \frac{\partial \varepsilon E}{\partial t},$$

$$\nabla \times E = -\frac{\partial \mu_e H}{\partial t},$$

$$p = p(\rho, T),$$

$$\frac{\partial \rho}{\partial t} + \frac{\partial}{\partial x^i}(\rho u^i) = 0,$$

$$\frac{\partial \rho u^i}{\partial t} + \frac{\partial}{\partial x^j}(\rho u^i u^j - \tau^{ij}) = F_e^i + F_g^i,$$

$$\frac{\partial \rho \overline{e_m}}{\partial t} + \frac{\partial}{\partial x^j}(p\overline{e_m} u^j - u^i \tau^{ij} - Q^j) = E^j \cdot J^j,$$

$$\frac{\partial \rho_e}{\partial t} + \frac{\partial J^j}{\partial x^j} = 0,$$

$$J^i = i^i + \rho_e u^i = \sigma[E^i + \mu_e(u \times H)^i] + \rho_e u^i.$$

## 3 Electromagnetic Equations in Plasma Dynamics

The fundament equations of plasma dynamics discussed in section 2 may be roughly divided into two groups; one deals with mainly the gasdynamics quantities, $T, p, \rho, u^i$ or their partial values, and the other deals with mainly the electromagnetic quantities, $E^i, H^i, \rho_e, J^i$, which may be called electromagnetic equations. Of course, it should be noted here that there are interaction terms in these two groups which are very important in our treatment and that electromagnetic equations of $\rho_e$ and $J^i$ are derived from the gasdynamic equations.

### 3.1 Generalized Ohm's Law.

As we saw in Section 2 the total electrical current $J$ consists of the conduction current $i$ and the convection current $\rho_e u$,

$$J = i + \rho_e u. \tag{1}$$

The equation of conservation of electrical charge of the plasma is given by (2.55) which may be written in the following form by the help of relation (2.58)

$$\frac{\partial \rho_e}{\partial t} + \frac{\partial}{\partial x^j}(\rho_e u^j) = -\frac{\partial i^j}{\partial x^j}. \tag{2}$$

It is interesting to notice the similarity of equation (3.2) with gasdynamic equation of continuity (2.27). The only difference occurs on the right hand side because in equation (2.27) we assume that there is no source or sink of mass for the plasma as a whole while in equation (3.2), the divergence of conduction current $i$ represents a source or sink of electric change in the plasma.

The equation (2.27) gives a differential equation for $J$, but instead of using it, we will find some simple relation of conduction current $i$, i.e.,

$$i = \sum^n \rho_{e\alpha} w_\alpha, \qquad (3)$$

in terms of some other quantities of plasma.

From equation (2.4) and (2.28) with the help of equation (2.1) and (2.27), we have the equation for the diffusion velocity $w_\alpha = u_\alpha - u$ as follow

$$\frac{\partial w_\alpha^i}{\partial t} + u_\alpha^j \frac{\partial u_\alpha^i}{\partial x^j} - u^j \frac{\partial u^i}{\partial x^j} = \frac{1}{\rho_\alpha} \cdot \frac{\partial \tau_\alpha^{ij}}{\partial x^j} - \frac{1}{\rho_\alpha} \cdot \frac{\partial \tau^{ij}}{\partial x^j} + \frac{X_\alpha^i}{\rho_\alpha} - \frac{X^i}{\rho} - \frac{\sigma_\alpha}{\rho_\alpha}(Z_\alpha^i - u_\alpha^i) \qquad (4)$$

If we want to find a relation for conduction current $i$ under the condition that it is explicitly independent of the time $t$ and spatial coordinates $x^j$, with the electromagnetic forces as the only predominant force and with no source term, i.e., $\sigma_\alpha = 0$, equation (3.4) reduces to the following simple form [Grad 1956],

$$\frac{F_{e\alpha}}{\rho_\alpha} = \frac{F_e}{\rho}, \qquad (5)$$

where

$$F_{e\alpha} = \rho_{e\alpha}[E + \mu_e(u_\alpha \times H)] + F_{e\alpha o} = \rho_{e\alpha}[E + (u_\alpha \times B)] + F_{e\alpha o},$$

$$F_e = \sum_{\alpha=1}^n F_{e\alpha} = \rho_e E + (J \times B) = \rho_e(E + (u \times B)) + (i \times B).$$

Thus from (3.5) we have

$$(\rho \rho_{e\alpha} - \rho_\alpha \rho_e)(E + u \times B) + (\rho \rho_{e\alpha} w_\alpha - \rho_\alpha i) \times B = -\rho F_{e\alpha o}, \qquad (6)$$

where we have used the relation $w_\alpha = u_\alpha - u$.

Usually we assume that the interaction force between charged particles of $\beta s$ with velocities $u_\alpha$ and $u_\beta$ is proportional to the difference of their velocities $u_\alpha - u_\beta$ with corresponding proportional constant $\alpha_{\alpha\beta}$, that is we assume that the interaction force $F_{e\alpha o}$ is given by

$$F_{e\alpha o} = \sum_{\beta=1}^n \alpha_{\alpha\beta}(u_\alpha - u_\beta) = \sum_{\beta=1}^n \alpha_{\alpha\beta}(w_\alpha - w_\beta). \qquad (7)$$

Using (3.7) in (3.6) gives

$$(\rho \rho_{e\alpha} - \rho_\alpha \rho_e)(E + u \times B) + (\rho \rho_{e\alpha} w_\alpha - \rho_\alpha i) \times B = -\rho \sum_{\beta=1}^n (w_\alpha - w_\beta), \qquad (8)$$

which is a linear equation for $w_\alpha$ because $i$ is a linear function of $w_\alpha$, as we have $i = \sum_{\alpha=1}^{n} \rho_{e\alpha} w_\alpha$. We may solve for $w_\alpha$ from equation (3.8).

For the simple case of fully ionized plasma equation (3.8) after neglecting small terms gives

$$i = \sigma(E + u \times B) - (\sigma\zeta)i \times B, \tag{9}$$

where $\sigma = -\rho_{e1}^2/\alpha_{12}$ is the electrical conductivity, and $\zeta = -1/\rho_{e1}$. In derivation of (3.9) we have assumed that $m_1 \ll m_2$ and $\rho_e \simeq 0$. The term $\sigma\zeta B$ is know as the Half vector and $\sigma\zeta i \times B$ is the Hall current. When the Hall current is negligible, equation (3.9) reduces to

$$i = \sigma(E + u \times B), \tag{9}$$

which has been widely accepted form of Ohm's law in the single fluid description of plasma. An interesting result from (3.9) is that

$$i \cdot (E + u \times B) = \frac{i^2}{\sigma}, \tag{10}$$

which is sometimes used as a definition for the electrical conductivity $\sigma$.

We would like to give a more precise form of Ohm's law for the case of a fully ionized plasma considering the quasi neutrality $\rho_e \simeq 0$. Thus let the plasma consist of two type of particles $\alpha = 1, 2$, and that

$$\rho_e = \rho_{e1} + \rho_{e2} = 0, \tag{11}$$

which gives

$$J = \rho_e u + i = i. \tag{12}$$

For simplicity in writing we use the parameters $\gamma_1, \gamma_2$ defined by

$$\gamma_\alpha = \frac{e_\alpha}{m_\alpha} = \frac{\nu_\alpha e_\alpha}{\gamma_\alpha m_\alpha} = \frac{\rho_{e\alpha}}{\rho_e}, \tag{13}$$

thus we have $\rho_{e\alpha} = \gamma_\alpha \rho_e$, and then equation (3.11) gives

$$\gamma_1 \rho_1 + \gamma_2 \rho_2 = 0. \tag{14}$$

Since we have $\rho = \rho_1 + \rho_2$ thus by (3.14) we have

$$\rho_1 = \frac{\gamma_2}{\gamma_2 - \gamma_1}\rho, \quad \rho_2 = \frac{\gamma_1}{\gamma_1 - \gamma_2}\rho. \tag{15}$$

Now the current is given by

$$i = J \cdot = \sum_{\alpha=1}^{2} \rho_{e\alpha} u_\alpha = \rho_{e1} u_1 + \rho_{e2} = \gamma_1 \rho_1 u_1 + \gamma_2 \rho_2 u_2$$
$$= (\gamma_1 \rho_1 + \gamma_2 \rho_2)u + \gamma_1 \rho_1 w_1 + \gamma_2 \rho_2 w_2$$
$$= (\gamma_1 \gamma_2/(\gamma_2 - \gamma_1))(w_1 - w_2)\rho. \tag{16}$$

or $w_1 = w_2 = \dfrac{\gamma_2 - \gamma_1}{\gamma_1 \gamma_2} \dfrac{i}{\rho}$. On the other hand we have

$$\rho u = \rho_1 + \rho_2 u_2 = \rho_1 u + \rho_2 u + \rho_1 w_1 + \rho_2 w_2, \tag{17}$$

thus

$$\rho_1 w_1 + \rho_2 w_2 = 0, \tag{18}$$

which by using (3.15) gives

$$\gamma_2 w_1 = \gamma_1 w_2. \tag{19}$$

Now solving (3.16) for $w_1$ or $w_2$ and using (3.19) gives

$$w_1 = -\dfrac{i}{\gamma_2 \rho}, \quad w_2 = -\dfrac{i}{\gamma_1 \rho}, \tag{20}$$

Hence

$$u_1 = u - \dfrac{i}{\gamma_2 \rho}, \quad U_2 = u - \dfrac{i}{\gamma_1 \rho}. \tag{21}$$

From (3.11), (3.13) and (3.15) we have

$$\rho_{e1} = \dfrac{\gamma_1 \gamma_2}{\gamma_2 - \gamma_1} \rho, \quad \rho_{e2} = \dfrac{\gamma_1 \gamma_2}{\gamma_1 - \gamma_2} \rho. \tag{22}$$

Using equations (3.13)-(3.21) in (3.8) gives

$$\rho \left( \dfrac{\gamma_1 \gamma_2}{\gamma_2 - \gamma_1} (E + u \times B) \right) - \dfrac{\gamma_1 + \gamma_2}{\gamma_2 - \gamma_1} i \times B = -\alpha_{12} \left( \dfrac{\gamma_2 - \gamma_1}{\gamma_1 \gamma_2} \right) \dfrac{i}{\rho},$$

or

$$i = -\left( \dfrac{\gamma_1 \gamma_2}{\gamma_2 \gamma_1} \right)^2 \dfrac{\rho^2}{\alpha_{12}} \cdot (E + u \times B) + \dfrac{(\gamma_1 + \gamma_2)\gamma_1 \gamma_2}{(\gamma_2 - \gamma_1)^2} \cdot \dfrac{\rho}{\alpha_{12}} (i \times B)$$
$$= \sigma(E + u \times B) - \chi(i \times B), \tag{23}$$

where $\sigma$, the electrical conductivity, and $\chi$, the Hall effect, are defined as follows

$$\sigma = -\left( \dfrac{\gamma_1 \gamma_2}{\gamma_2 - \gamma_1} \right) \cdot \dfrac{\rho^2}{\alpha_{12}}, \chi = \dfrac{(\gamma_1 + \gamma_2)\gamma_1 \gamma_2}{(\gamma_2 - \gamma_1)^2} \cdot \dfrac{\rho^2}{\alpha_{12}}, \tag{24}$$

where $\alpha_{12}$ is the friction coefficient and should be calculated from kinetic theory of plasma. If we assume that the interaction force is simply due to the momentum gained by electrons in collision with ions, the simple kinetic theory gives

$$F_{e10} = m_1 \nu_1 f_1 (u_2 - u_1) = \alpha_{12}(u_1 - u_2) = -F_{e20} \tag{25}$$

Hence

$$\alpha_{12} = -m_1 \nu_1 f_1 = -m_2 \nu_2 f_2 = \alpha_{21} \tag{26}$$

where $f_1$ is the collision frequency of each electron. We have then

$$\sigma = \left(\frac{\gamma_1 \gamma_2}{\gamma_2 - \gamma_1}\right)^2 \frac{\rho^2}{m_1 \nu_1 f_1}, \chi = \frac{(\gamma_1 + \gamma_2)\gamma_1 \gamma_2}{(\gamma_2 - \gamma_1)^2} \cdot \frac{\rho}{m_1 \nu_1 f_1}. \quad (27)$$

If we assume the plasma is singly ionized, i.e.,

$$e_1 = -e_2 = e \quad (28)$$

then from (3.11) we obtain $\nu_1 = \nu_2 = \nu$ and that $\gamma_1 = \dfrac{e}{m_1}, \gamma_2 = \dfrac{-e}{m_2}$. $\sigma$ and $\chi$ can be expressed in terms of $m_1, m_2$ as follows using these,

$$\sigma = \frac{e^2 \nu}{m_1 f_1} = e^2 \cdot \frac{\nu}{m_1 f_1}, \chi = e \frac{m_1 - m_2}{m_1 + m_2} \cdot \frac{1}{m_1 f_1}, \quad (29)$$

where we have used $\rho = (m_1 + m_2)\nu$. If we assume that the collision frequency of each electron is proportional to the number density of electrons, say $f_1 = c_1 \nu_1 = c_1 \nu$, then $\sigma$ is constant and is given by

$$\sigma = e^2 \frac{1}{m_1 c_1} = \text{constant}. \quad (30)$$

If the plasma is partial ionized equation (3.8) can be used to obtain expressions for $i$ in terms of $E$, $u$ and $B$, see [Grad 1956] or [Pai 1963].

## 3.2   Lorentz Force.

Most of the equation describing the motion of plasma given in section 2 cure in the form of conservation laws, except the terms $\sigma_\alpha$ and $X_\alpha^i$ which are not represented in the form of divergence of a quantity. In the next section for deriving the shock layer equations it will be very useful to express these equations in the form of conservation laws. We assume that $\sigma_\alpha = 0$, $X_\alpha = F_{e\alpha}$ and $F_{g\alpha} = 0$. The terms $F_e = \rho_e E + (J \times B)$ and E.J on the rhs of equations (2.28) and (2.4) are not in the form of conservation laws, but can be expressed in that form as follows. By Maxwell's equations we have

$$\rho_e E + J \times B = (\nabla \cdot E)E + (\nabla \times \frac{B}{\mu_e}) \times B - \frac{\partial \varepsilon E}{\partial t} \times B, \quad (31)$$

then writing

$$B \times \frac{\partial \varepsilon E}{\partial t} = -\frac{\partial}{\partial t}(E \times B) + E \times \frac{\partial B}{\partial t}$$

and adding $B(\nabla \cdot B) = 0$ to the rhs of equation (3.31), we obtain

$$\rho_e E + (j \times B) = [E(\nabla \cdot E) + \frac{B}{\mu_e}(\nabla \cdot B) - E \times (\nabla \times E) - B \times (\nabla \times \frac{B}{\mu_e})] - \frac{\partial}{\partial t}(E \times B),$$

which after some calculations gives

$$\rho_e E^i + (J \times B)^i = \frac{\partial}{\partial x^j}\left[\varepsilon E^i E^j + \frac{1}{\mu_e}B^i B^j - \frac{1}{2}(\varepsilon E \cdot E + \frac{1}{\mu_e}B \cdot B)\delta^{ij}\right] - \frac{\partial}{\partial t}(\varepsilon E \times B)^i. \quad (32)$$

with a similar calculation we obtain the following formula for $E.J$,

$$\begin{aligned} E.J &= -\frac{\partial}{\partial \chi^j}(E \times \frac{B}{\mu})^j - \frac{\partial}{\partial t}(\frac{1}{2}\varepsilon E.E + \frac{1}{2\mu_e}B.B) \\ &= -\text{div}(E \times \frac{B}{\mu}) - \frac{\partial}{\partial t}(\frac{1}{2}\varepsilon E.E + \frac{1}{2\mu_e}B.B). \end{aligned} \quad (33)$$

## 4   Conservation laws and shock wave structure

Under standard physical assumptions, the continuum equations governing the motion of plasma in space $\chi \in \mathbb{R}$ and $t \in \mathbb{R}^+$ is in the form of a system of parabolic pdes which can be stated in the following general form

$$u_t(x,t) + f(u(x,t))_x = (A(u,\lambda)u_x)_x, u \in \mathbb{R}^n \quad (1)$$

where the flux $f$ and the viscosity matrix $A$ are smooth, $\lambda \in [0,\infty)^m$ is the viscosity vector, comprising of dissipation coefficient, and $A(u,0) = 0$. In the limiting case $\lambda = 0$, (4.1) is the system of conservation laws

$$u_t(x,t) + f(u(x,t))_x = 0, \quad (2)$$

describing ideal plasma with no dissipation effect taken into account. The conservation laws (4.2) admit certain discontinuous weak solutions of the form

$$u(u^+, u^-, s) := u(x,t) = \begin{cases} u^-, & x < st, \\ u^+, & x > st, \end{cases} \quad (3)$$

which are called shock waves, with speed $s$, and satisfy the Rankin-Hugoniot jump conditions

$$f(u^+) - su^+ = f(u^-) - su^- = C, \quad (4)$$

where $C$ is a constant. As a result from (4.4) we have

$$[f(u)] := f(u^+) - f(u^-) = s(u^+ - u^-) = s[u]$$

and since $u^+, u^-$ are constant we see that $s$, the propagation speed of the shock wave must be constant,

$$s = \text{constant}. \quad (5)$$

The values $u^-$ and $u^+$ correspond to the two states of the plasma connected by the abrupt discontinuity given by equation (4.3). Usually $u^-$ is called the upstream state, in behind, of the shock and $u^+$ is called the downstream state, in front, of the shock. The interface of the shock wave is $x = st$ which moves with velocity $s$ in the $x$-space. Principal for our work is the notion of structure of viscouse profile for a shock wave. There two approach to this notion, one from physical point of view and the other from mathematical point of view. Physically it is very essential to specify which type of shock wave will occur in nature, and which type will not occur. Physical a shock wave is a rapid change occurring in a thin layer which propagates in the medium. Thus a shock wave to be a thin transition region across which the macroscopic variables leap from their upstream values to their downstream values, we shall take the shock front to be planar, an orthogonal to the $x$-axis, and choose the flow direction such that $x = -\infty$ is upstream and $x = +\infty$ is downstream. Strictly, the shock extends from one limit to the other, but on macroscopic scale almost all the change occurs in a thin region of thickness $\lambda_s$. In the following we shall take the shock thickness $\lambda_s$ to be negligible, and deal only with the relations between the variables upstream of the shock and those downstream of it. If a shock wave merits to occur naturally then in fact it is a limit of the solutions $u(x,t)$ for the pde (4.1). Moving in a frame attached to the shock interface we see that $u(x,t)$, the vector of variables of the plasma, must depend only on $\zeta := x - st$ as the shock wave (4.3) is. Hence we have to see if there is a solution $u(x,t)$ of (4.1) in the form $u(\zeta) = u(x - st) = u(x,t)$ which tends to $u^-$ as $x \to -\infty$ and tends to $u^+$ as $x \to +\infty$. It can be easily seen that this is equivalent to ask if there is a solution $u(\zeta) = u(x - st) = u(x,t)$ of (4.1) which tends to $u^+$ and $u^-$ as $\zeta \to -\infty$ and $\zeta \to +\infty$, respectively. Physically this means that in a reference frame attached to the shock interface, and hence moving with the constant velocity $s$, we have $s = 0$, so $x = \zeta$ in this reference frame. Motivated by this explanations we define the structure or viscouse profile for the shock wave (4.3) as follows.

**Definition 4.1:** A structure or viscouse profile for the shock wave

$$u(u, u^+, s) = \begin{cases} u^-, & x < st \\ u^+, & x > st, \end{cases}$$

is a travelling wave solution $u(x,t) = u(x - st)$ which tends to the limits $u^-, u^+$ as $x \to -\infty$, $x \to -\infty$, respectively.

Let's to examine more closely what will happen if there is a structure $u(\zeta) = u(x - st)$ for the shock wave (4.3). First of all $u(\zeta) = u(x - st)$ must be a solution of pde (4.1), since we have $\partial/\partial x = d/d\zeta$ and $\partial/\partial t = -sd/d\zeta$, this gives

$$-s\frac{du}{d\zeta} + \frac{df(u)}{d\zeta} = \frac{d}{d\zeta}\left(B(u,\lambda)\frac{du}{d\zeta}\right),$$

and integrating this, we obtain

$$-su + f(u) + C = B(u,\lambda)\frac{du}{d\zeta},$$

which is a system of odes that the structure $u(\zeta) = u(x - st)$ must satisfy. Here $C$ is the constant of integration. Now considering the two boundary conditions $u(\zeta) \to u^{\pm}$ as $\zeta \to \pm\infty$, and the fact that $u(\zeta)$ is a solution of the above system of odes we obtain $du(\zeta)/d\zeta = 0$ at $\zeta = \pm\infty$. Thus we must have

$$C = su^- f(u^-) = su^+ - f(u^+). \tag{6}$$

Using (4.6) we see that $u^-$ and $u^+$ are both rest points of the above system of odes. In summary we see that the structure $u(\zeta) = u(x - st)$ is a solution of the following system of odes

$$B(u,\lambda)\frac{du}{d\zeta} = G(u), \tag{7}$$

$$G(u) = f(u) - su + C, \ G(u^-) = G(u^+) = 0,$$

which connects the rest point $u^- = u(-\infty)$ to the rest point $u^+ = u(+\infty)$. The system of odes (4.7) is called the shock layer equations. A solution of the odes connecting two different rest points are called heteroclinic orbits. We state the above arguments in the form of a theorem as follows,

**Definition 4.2:** A structure for the shock wave

$$u(u^-, u^+, s) = \begin{cases} u^-, & x < st, \\ u^+ & x > st, \end{cases}$$

is a travelling wave solution $u(x,t) = u(x - st) - u(\zeta)$ of the parabolic system of odes

$$u_t(x,t) + f(u(x,t))_x = (B(u,\lambda)u_x)_x,$$

which tends to the limit $u^-, u^+$ as $\zeta \to -\infty$, $\zeta \to +\infty$, respectively, must be a heteroclinic orbit of system of shock layer equations

$$B(u,\lambda)\frac{du}{d\zeta} = f(u) - su + C$$

running from the upstream rest point $u^- = u(-\infty)$ to the downstream rest point $u^+ = u(+\infty)$.
□

Some authors use this theorem to define the structure of a shock wave.

Now we look, at the problem of specifying shock wave solutions (4.3) for the system of conservation laws (4.2), from a different point of view. Mathematically shock waves are discontinuous weak solutions of the system of conservation laws (4.2). The disadvantage of the shock waves as weak solutions for (4.2) is that they suffer from non uniqueness, there may be even uncountably many shock wave weak solutions for equation (4.2) with a nice boundary condition. The question naturally arises which criteria one should apply to specify proper solutions of the conservation laws. Despite there some criteria such as Olenik's entropy condition and Lax's entropy condition, most widely accepted is the viscouse profile criterion. A shock wave solution (4.3) satisfies this criterion if there exists, inevitably, travelling wave solution $u_\varepsilon(x,t) = u_\varepsilon(x - st/\varepsilon)$ of the system of conservation laws with viscosity term added

$$u_t(x,t) + f(u(x,t))_x = \varepsilon(D(u,\mu)u_x), \tag{8}$$

which converges to $u^-$ for $x - st < 0$ and converges to $u^+$ for $x - st > 0$ as $\varepsilon \to 0$. We state this in the form of a definition.

**Definition 4.3:** A shock wave solution (4.3) is a viscosity solution or admit a $D(u,\mu)$-viscouse profile if there exist travelling wave solutions $u_\varepsilon(x - st/\varepsilon)$ of (4.8) which converge to the shock wave (4.3) as $\varepsilon \to 0$. The travelling wave solution $u_\varepsilon(x - st/eps)$ is called a viscouse profile for the shock wave (4.3).

Let $u_\varepsilon(x - st\varepsilon)$ be a viscouse profile for the shock wave (4.3), then simple computations and setting $u(\eta) = u_\varepsilon(x - st/\varepsilon)$ shows that $u(\eta)$ is heteroclinic orbit of the system of odes

$$D(u,\eta)\frac{du}{d\eta} = f(u) - su - (f(u^-) - su^-) = G(u),$$
$$G(u^-) = G(u^+) = 0. \tag{9}$$

Similar to Theorem 4.1 we have the following theorem.

**Theorem 4.2:** A viscouse profile for the shock wave (4.3) is a heteroclinic orbit $u(\eta)$, at the system (4.9), connecting the rest point $u^- = u(-\infty)$ to the rest point $u^+ = u(+\infty)$. □

By considering Theorems 4.1, 4.2 we see that both structure and viscouse profile for the shock wave (4.3) are heteroclinic orbits of the system of shock layer equations (4.7) and (4.9), respectively. Thus from now on we identify the heteroclinic orbit of (4.7) or (4.9) running from $u^-$ to $u^+$ with structure or viscouse profile for shock wave (4.3) and use the terms structure and viscouse profile as synonyms. Here we should mention that the two system of parabolic pdes (4.1) and (4.8) and the corresponding notions structure and viscouse profile are identical can be seen by changing $x, t$ in (4.1) to $\varepsilon^{-1}x, \varepsilon^{-1}t$ gives the system (4.8) and accordingly the definition of structure reduces to that of the viscouse profile. As a consequence of the above arguments and theorems we can state the following equivalent definition for structure or viscouse profile for the shock wave (4.3).

**Equivalent Definition:** A structure or viscouse profile for the shock wave (4.3) is a heteroclinic orbit of the system of shock layer equations running from the rest point $u^-$ to the rest point $u^+$

Now that we have clearly declared the notion of structure for a shock wave, we can impose more extra conditions on the structure for a shock wave to select, more distinctively, among the shock wave solutions of the conservation laws. Looking at the geometry of the structure of shock waves, as defined in Definition 4.1, we see that the convergence of the travelling wave solution to the shock wave is essentially a condition at infinity, and at the shock wave interface it may occur that the travelling wave solution differ greatly from the shock wave. Here note that the travelling wave solution $u(\zeta)$ of (4.1) depends on the viscosity vector $\lambda$ which describes the dissipations.

Let's denote the structure $u(\zeta)$ which is a heteroclinic orbit of the shock layer equations (4.7) corresponding to $\lambda$ by $u(\zeta, \lambda)$. Suppose that the heteroclinic orbit $u(\zeta, \lambda)$ exits for sufficiently small $|\lambda|$ and that $u(\zeta, \lambda)$ converges uniformly to the rest points $u^-, u^+$ for $\zeta < 0, \zeta > 0$ as $|\lambda| \to 0$. We state this in the form of a definition as follow. The reader should beware not to confuse the following definition of stability of a structure or viscouse profile with many other similar notions, see eg. [ ].

**Definition 4.3:** A structure $u(\zeta)$ for the shock wave (4.3) is stable (in the sense of Germain) if the heteroclinic orbit $u(\zeta)$ exits for all sufficiently small $|\lambda|$ and that $u(\zeta)$ converges uniformly to the shock wave (4.3) as $|\lambda| \to 0$.

## 4.1 Isotropy, Gallilean invariance and Rotational invariance.

The laws of continuum mechanics are subject to the Gallilean invariance which is the rule of classical mechanics. Consider a frame $\mathcal{R}'$ moving with uniform speed $V$ w.r.t. the reference frame $\mathcal{R}$; thus the new independent variables in $\mathcal{R}'$ are

$$x' = x - vt, \quad t' = t.$$

Looking at the dependent variables in the frame $\mathcal{R}'$, the density $\rho_\alpha$, the pressure $p_\alpha$, the magnetic intensity $B$, the heat flux $Q_\alpha$ and the internal energy $e_\alpha$ are invariant,

$$\rho'_\alpha = \rho_\alpha, \ e_\alpha = e_\alpha, \ p_\alpha = p_\alpha,$$

while the velocity $u_\alpha$, the electric field $E$ and the current density $J$ are transformed as follows,

$$u'_\alpha = u_\alpha - v, \ E' = E + v \times B, J' = J + \rho_e v.$$

The expression Gallilean invariance means that the equations describing the motion of the plasma in the reference frame $\mathcal{R}$ remain valid in the preference frame $\mathcal{R}'$ with corresponding dependent variables. The invariance of the equations is easily established using the rules

$$\frac{\partial}{\partial x^j} = \frac{\partial}{\partial x'^j}, \ \frac{\partial}{\partial t'} + v.\nabla_x = \frac{\partial}{\partial t} + v^j \frac{\partial}{\partial x^j}$$

Let us check, for instance, the equation of continuity of $\alpha_s$:

$$\begin{aligned}\frac{\partial \rho_\alpha}{\partial t'} + \frac{\partial \rho_\alpha u'^j_\alpha}{\partial x'^j} &= \frac{\partial \rho_\alpha}{\partial t} + v^j.\frac{\partial \rho_\alpha}{\partial x^j} + \frac{\partial(\rho_\alpha(u^j_\alpha - v^j))}{\partial x^j} \\ &= \frac{\partial \rho_\alpha}{\partial t'} + v^j.\frac{\partial \rho_\alpha}{\partial x^j} + \frac{\partial \rho_\alpha u^j_\alpha}{\partial x^j} - v^j \frac{\partial \rho_\alpha}{\partial x^j} \\ &= \frac{\partial \rho_\alpha}{\partial t} + \frac{\partial \rho_\alpha u^j_\alpha}{\partial x^j}.\end{aligned}$$

The computations are similar.

As we mentioned before, the Rankin-Hoguniot jump conditions for weak shock wave solutions of conservation laws implies that $s$, the speed of shock, is constant. Let the frame $\mathcal{R}'$ move with velocity $sn$ w.r.t. the reference frame $\mathcal{R}'$, where $n$ is the unit normal to the surface of the shock. Now we write the equations of motion of the plasma in the frame $\mathcal{R}'$. Then in the frame $\mathcal{R}'$ the shock surface has velocity zero, thus the equation of the shock layer in this case will be reduced to the following form

$$B(u,\lambda)\frac{du}{d\zeta} = f(u) + C. \tag{10}$$

This form for the shock layer equations can also be obtained by assuming a steady state motion for plasma, that is assuming $\frac{\partial}{\partial t} = 0$ and that the variables depend only on one space coordinate, say $x_0$. With this assumption the pde (4.1) reduces to

$$\frac{d}{dx}f(u) = \frac{d}{dx}(B(u,\lambda)\frac{du}{dx})$$

which by integration gives

$$B(u,\lambda)\frac{du}{dx} = f(u) + C \tag{11}$$

which is equivalent to equation (4.10). This simple approach for deriving the shock layer equations is adapted by many authors and has the advantage that it restores the coordinate variable $x$, and hence make more physical insight to the shock layer equations.

The next simplifying fact which should be mentioned here is the rotational invariance of the equation of motion of plasma. It is known that the velocity vector and the magnetic field intensity are invariant under a spatial rotation. This fact is important in the implementations for numerical algorithms in conservation laws as well as in the study of stability of some types of shock waves.

We will not go further in detail of the analysis of rotational invariance, the interested reader can consult Freistuhler [ , ] and the book [ ] where an explicit computation is given. The implication of the rotational invariance in our work is that we can assume that the direction of motion of the shock wave, i.e., the normal to the surface of the shock, is parallel to the $x$-axis. Then using the Galilean invariance we obtain the shock layer equations as given in (4.11). Finally there is another possibility for simplifying the equations, which is simply to take units of measurement so that some constant parameters have suitable values, e.g., we may take $\mu_e = \varepsilon = 1$. This values will be used in section.

# 5  Shock layer equations

In the previous section we saw that the shock layer equations can be obtained by taking

$$\frac{\partial}{\partial t} = \frac{\partial}{\partial y} = \frac{\partial}{\partial z} = 0, \ x = x', y = x^2, z = x^3 \tag{1}$$

in the partial differential equations of motion of the plasma. This give the shock layer equations in the form of equation (4.11). Our aim is to obtain the shock layer equations for a fully ionized plasma, and we assume that there is no chemical reactions or ionization. Hence we assume that there is only two type of particles the ions with mass $m_i$ and the charge $e_i$ and the electrons with mass $m_e$ and the charge $e_e$. Let us denote $x'$ by the continuity equation (2.2) gives

$$\frac{d}{dx}(\rho_\alpha u_\alpha^1) = 0, \ \alpha = i, e \tag{2}$$

where we have used (5.1) and $\sigma_\alpha = 0$. Similarly from equation (2.4) we get

$$\frac{d}{dx}(\rho_\alpha u_\alpha^i u_\alpha^1 - \tau_\alpha^{i1}) = X_\alpha^i, \ \alpha = i, e \tag{3}$$

where $X_\alpha^i$ is given by (2.b). We assume that there is no gravitational forces, $F_{g\alpha} = 0$. Moreover since there is only two types of particles, hence using (2.8) and (3.4) we have

$$F_{eo} = -F_{eio} = \alpha_{12}(u_e - u_i) = -\alpha_{12}(u_i - u_e). \tag{4}$$

Therefore by (2.6) and (2.7) we see that

$$\begin{aligned} X_e &= \rho_e e[E + (u_e \times B)] + \alpha_{12}(u_e - u_i), \\ X_i &= \rho_e i[E + (u_i \times B)] - \alpha_{12})u_e - u_i). \end{aligned} \tag{5}$$

The equation of conservation of energy (2.9) for each species gives

$$\frac{d}{dx}(\bar{e}_\alpha u_\alpha^1 - u_\alpha^i \tau_\alpha^{i1} - Q_\alpha^1) = \varepsilon_\alpha, \tag{6}$$

where by (2.11a), $\varepsilon_{c\alpha} = 0$ and (2.11b) we have

$$\varepsilon_\alpha = \rho_{e\alpha} u_\alpha . E = J_\alpha . E, \tag{7}$$

and by (2.11) we have

$$Q_\alpha = Q_{c\alpha} + Q_{R\alpha}. \tag{8}$$

The Maxwell equations (2.12), (2.13) using (5.1) give,

$$\frac{dH_3}{dx} = -J_2, \quad \frac{dH_2}{dx} = J_3, \quad O = J_1,$$
$$\frac{dE_2}{dx} = \frac{dE_3}{dx} = 0. \tag{9}$$

By virtue of $\nabla . E = 0$ and $\nabla . B = 0$, and using (5.1) we obtain

$$\frac{dE'}{dx} = 0 \quad \frac{dB^1}{dx} = 0. \tag{10}$$

Note that the equations and results of Subsection 3.1 can be used here, where we may assume that the subscript 1 denote the ions and the subscript 2 denote that electrons. Thus from (3.12), (3.21) and (5.9) we have

$$J = i, \quad i' = 0, \quad u_e^1 = u_i^1 = u^1. \tag{10}$$

Equations (5.2), (5.3), (5.5) and (3.11) gives

$$\frac{d}{dx}(\rho_\alpha u_\alpha^1 u_\alpha^i - \tau_\alpha^{i1}) = \frac{d}{dx}(\rho_\alpha u_\alpha^1 + \rho_\alpha u_\alpha^1 \frac{du_\alpha^i}{dx} - \frac{d\tau_\alpha^{i1}}{dx}$$
$$= \rho_\alpha u^1 \frac{du_\alpha^i}{dx} - \frac{d\tau_\alpha^{i1}}{dx}$$

therefore

$$\rho_e u^1 \frac{du_e^i}{dx} - \frac{d\tau_e^{i1}}{dx} = \rho_{ee}[E^i + (u_e \times B)^i] + \alpha_{12}(u_e^i - u_i^i),$$
$$\rho_i u^1 \frac{du_i^i}{dx} - \frac{d\tau_i^{i1}}{dx} = \rho_{ei}[E^i + (u_i \times B)^i] - \alpha_{12}(u_e^i - u_i^i), \tag{11}$$

Thus we have

$$\rho_e u^1 \frac{du_e^i}{dx} + \rho_i u^1 \frac{du_i^1}{dx} - \frac{d(\tau_e^{i1} + \tau_i^{i1})}{dx} = [(\rho_{ee} u_e + \rho_{ei} u_i) \times B]^i$$
$$= [i \times B]^i, \tag{12}$$

where we have used (5.1), (3.11) and (3.12). By the above results we have $i \times B = (i_2 B_3 - i_3 B_2, i_3 B_1, -i_2 B_2)$. Taking appropriate linear combinations of (5.11), we have by ( . )

$$\begin{aligned}
\frac{d}{dx}(J^i u^1) &= \frac{d}{dx}(\gamma_e \tau_e^{i1} + \gamma_i \tau_i^{i1}) + (\gamma_e \rho_{ee} + \gamma_i \rho_{ei}) + E^1 + [\gamma_e \rho_{ee} u_e + \gamma_i \rho_{ei} u_i) \times B]^i \\
&\quad + (\gamma_e - \gamma_i) \alpha_{12}(u_e^i - u_i^i) \\
&= \frac{d}{dx}(\gamma_e \tau_e^{i1} + \gamma_i \tau_i^{i1}) + \gamma_e \gamma_i \rho E^i - \gamma_e \gamma_i \rho (u \times B)^i + (\gamma_e + \gamma_i)(J \times B)^i \\
&\quad + (\gamma_e - \gamma_i) \alpha_{12}(u_e^i - u_i^i).
\end{aligned} \qquad (13)$$

This equation for $i = 1, 2, 3$ gives

$$0 = \frac{d}{dx}(\gamma_e \tau_e^{11} + \gamma_i \tau_i^{11}) - \gamma_e \gamma_i \rho E^1 - \gamma_e \gamma_i \rho (u^2 B^3 - u^3 B^2) - \frac{1}{2}(\gamma_e + \gamma_i) \frac{d}{dx}((B^2)^2 + (B^3)^2),$$

$$\frac{d}{dx}\left(-u^1 \frac{dB^3}{dx}\right) = \frac{d}{dx}(\gamma_e \tau_e^{21} + \gamma_i \tau_i^{21}) - \gamma_e \gamma_i \rho E^2 - \gamma_e \gamma_i \rho (u^1 B^3 - u^3 B^1) + (\gamma_e + \gamma_i) B^1 \frac{dB^2}{dx}$$
$$+ \frac{\alpha_{12}}{\rho} \frac{(\gamma_e - \gamma_i)^2}{\gamma_e \gamma_i} \frac{dB^3}{dx},$$

$$\frac{d}{dx}\left(u^1 \frac{dB^2}{dx}\right) = \frac{d}{dx}(\gamma_e \tau_e^{31} + \gamma_i \tau_i^{31} + \gamma_i \tau_i^{31}) - \gamma_e \gamma_i \rho E^3 - \gamma_e \gamma_i \rho (u^1 B^2 - u^2 B^1)$$
$$- (\gamma_e + \gamma_i) B^1 \frac{dB^3}{dx} - \frac{\alpha_{12}}{\rho} \frac{(\gamma_i - \gamma_i)^2}{\gamma_e \gamma_i} \frac{dB^2}{dx}.$$

Multiplying by $(\gamma_e \gamma_i \rho)^{-1}$ gives

$$0 = \frac{1}{\gamma_\alpha \gamma_i \rho} \frac{d}{dx}(\gamma_e \tau_e^{11} + \gamma_i \tau_i^{11}) - E^1 - (u^2 B^2 - u^3 B^2) - \frac{1}{2} \frac{\gamma_e + \gamma_i}{\gamma_e \gamma_i \rho} \frac{d}{dx}((B^2)^2 + (B^3)^2)$$

$$\frac{1}{\gamma_e \gamma_i \rho} \frac{d}{dx}\left(-u^1 \frac{dB^3}{dx}\right) = \frac{1}{\gamma_e \gamma_i \rho} \frac{d}{dx}(\gamma_e \tau_e^{21} + \gamma_i \tau_i^{21}) - E^2 - (u^1 B^3 - u^3 B^1)$$
$$+ \frac{\gamma_e + \gamma_i}{\gamma_e \gamma_i \rho} B^1 \frac{dB^2}{dx} + \frac{\alpha_{12}}{\rho^2} \frac{(\gamma_e - \gamma_i)^2}{(\gamma_e \gamma_i)^2} \frac{dB^3}{dx} \qquad (13)$$

$$\frac{1}{\gamma_e \gamma_i \rho} \frac{d}{dx}\left(u^1 \frac{dB^2}{dx}\right) = \frac{1}{\gamma_e \gamma_i \rho} \frac{d}{dx}(\gamma_e \tau_e^{31} + \gamma_i \tau_i^{31}) - E^3 - (u^1 B^2 - u^2 B^1)$$
$$- \frac{\gamma_e - \gamma_i}{\gamma_e \gamma_i \rho} B^1 \frac{dB^3}{dx} - \frac{\alpha_{12}}{\rho^2} \frac{(\gamma_e - \gamma_i)^2}{(\gamma_e \gamma_i)^2} \frac{dB^2}{dx}.$$

Now we turn to the equation of energy for each species, From equation (5.6) and (5.7) we have, for $\alpha = i, e$

$$\frac{d}{dx}\left(\rho_\alpha u_\alpha^1 (U_{m\alpha} + \frac{1}{2} u_\alpha^2 + \phi_\alpha + \frac{E_{R\alpha}}{\rho_\alpha}) - u_\alpha^i \tau_\alpha^{i1} - Q_\alpha^1\right) = \rho_{e\alpha} u_\alpha . E \qquad (14)$$

which using (5. ), (5. ) and (5. ) gives

$$\frac{d}{dx}\left(\rho u^1 \frac{\gamma_e}{\gamma_e - \gamma_i}(U_{mi} + \frac{1}{2}(u)^2 + \frac{1}{2}\frac{J^2}{\gamma_e^2 \rho^2} - \frac{u.J}{\gamma_e \rho} + Q_e) + E_{Re} - u^i \tau_i^{i1} + \frac{1}{\gamma_i \rho} J^i \tau_i^{i1} - Q_e^1\right) \rho_u = \frac{\gamma_e \gamma_i \rho}{\gamma_i - \gamma_e} u.E \qquad (15)$$

$$\frac{d}{dx}\left(\rho u^1 \frac{\gamma_e}{\gamma_e - \gamma_i}(U_{me} + \frac{1}{2}(u)^2 + \frac{1}{2}\frac{J^2}{\gamma_e^2 \rho^2} - \frac{u.J}{\gamma_e \rho} + Q_i)\right.$$
$$\left. + E_{Ri} - u^i \tau_i^{i1} + \frac{1}{\gamma_e \rho} J^i \tau_i^{i1} - Q_i^1 \right) \tag{16}$$
$$= \frac{\gamma_i \gamma_i \rho}{\gamma_e - \gamma_i} u.E - \frac{\gamma_i}{\gamma_e - \gamma_i} J.E,$$

Adding these two equations we get

$$\frac{d}{dx}(\rho u^1 \left(\frac{\gamma_i(U_{me} + \phi_e) - \gamma_e(U_{mi} + \phi_i)}{\gamma_i - \gamma_e}\right) - \frac{J^2}{2\gamma_e \gamma_i \rho^2} + \frac{1}{2}\rho u^1(u)^2 + E_{Re} - u^i(\tau_e^{i1} + \tau_i^{i1})$$
$$- \frac{J^i}{rho}(\frac{\tau_e^{i1}}{\gamma_i} + \frac{\tau_i^{i1}}{\gamma_e}) - (Q_e^1 + Q_i^1)) = J.E = J^2 E^2 + J^3 E^3$$
$$= -E^2 \frac{dB^3}{dx} + E^3 \frac{dB^3}{dx}.$$

Integrating this and using (5. ) gives

$$\rho u^1 \left(\frac{\gamma_i(U_{me} + \phi_e) - \gamma_e(U_{mi} + Q_i)}{\gamma_i - \gamma_e}\right) - \frac{J^2}{2\gamma - e\gamma_i \rho^2}$$
$$+ \frac{1}{2}\rho u^1(u)^2 + E_{Re} + E_{Ri} - u^i(\tau_e^{i1} + \tau_i^{i1})$$
$$- \frac{J^i}{\rho}\left(\frac{\tau_e^{i1}}{\gamma_e} + \frac{\tau_i^{i1}}{\gamma_e}\right) - (Q_e^1 + Q_i^1) = E^3 B^2 - E^2 B^3 + C, \tag{17}$$

where $C$ is a constant. It may seem that the first term in (5.17) is some what irrelevant, but it can be easily seen that it is the mean total internal energy of the plasma as a whole, as follows

$$\frac{1}{\gamma_i - \gamma_e}(\gamma_i(U_{me} + \phi_e) - \gamma_e(U_{mi} + \phi_i)) = \frac{m_e m_i}{m_e + m_i}\left(\frac{1}{m_i}(U_{me} + \phi_e) + \frac{1}{m_e}(U_{mi} + \phi_i)\right)$$
$$= \frac{m_e(U_{me} + \phi_e) + m_i(U_{mi} + \phi_i)}{m_e + m_i} = U_m + \phi.$$

Therefore equation of conservation of energy, (5.17), can be written in the following form

$$\rho u^1(U_m + \phi) + E_{Re} + E_{Ri} - \frac{J^2}{2\gamma_i \gamma_e \rho^2} + \frac{1}{2}\rho u^1(u)^2 - u^i \tau^{i1}$$
$$- \frac{J^i}{\rho}\left(\frac{\tau_e^{i1}}{\gamma_i} + \frac{\tau_i^{i1}}{\gamma_e}\right) - (Q_e^1 + Q_i^1) \tag{18}$$
$$= E^3 B^2 - E^2 B^3 + C.$$

This equation can be obtained from equation (2.40), as follows

$$\frac{d}{dx}(\rho u^1(U'_m + \phi + \frac{1}{2}(u)^2) + E_R - u^i \tau^i - Q) = E.J, \tag{19}$$

where we have used (2.49) and (2.50), and $U'_m$ is give by (2.41a).

Then using (2.41), (3.15) and (3.20) we have

$$\begin{aligned} U'_m &= U_m + \frac{1}{2}\left(\frac{\gamma_i}{\gamma_i - \gamma_e}w_e^2 + \frac{\gamma_e}{\gamma_e - \gamma_i}w_i^2\right) \\ &= U_m + \frac{1}{2(\gamma_i - \gamma_e)}(\gamma_i w_e^2 - \gamma_e w_i^2) \\ &= U_m + \frac{1}{2(\gamma_i - \gamma_e)}\left(\gamma_i(\frac{J}{\gamma_e \rho}) - \gamma_e \frac{J}{\gamma_e \rho})^2\right) \\ &= U_m + \frac{1}{2(\gamma_i - \gamma_e)}\left(\frac{1}{\gamma_i} - \frac{1}{\gamma_e}\right)\frac{J^2}{\rho^2} \\ &= U_m - \frac{1}{2\gamma_e \gamma_i \rho^2}J^2. \end{aligned} \qquad (20)$$

Thus from (5.19) we get

$$\rho u^1(U_m + \phi) + \rho u^1(u)^2 + E_R - u^i \tau^{i1} - \frac{\rho u^1}{2\gamma_e \gamma_i \rho^2}J^2 - Q^1 = E^3 B^2 - E^2 B^3 + C.$$

At this point only a few steps remain to reach at our stated goal in this section, but however some remarks are in order. Looking at all equations which we have obtained so far we see that there is essentially two way to obtain the shock layer equations for the two fluid model for the fully ionized plasma. The first approach is simply writing the equations in term of partial variables. In this approach we have to work with many equations which depend on electron fluid and ion fluid dissipation coefficients, such as kinematic and transverse viscosities and the heat conductivity of each fluid. The second and more preferred approach is to work with mean variables of the plasma and then using equations (3.11)-(3.30) obtain the partial variables in terms of the mean variables. But there is one disadvantage in this approach, namely there is no explicit relation between $T_i, T_e$ and $T$, in such a way that when $T$ is determined we can obtain $T_i$ and $T_e$. Fortunately there are relations between $T_i, T_e$ and $T$ which we have not mentioned so far, that is the Newton's cooling law. By this law we have the following differential equation for $T_i$, [Shkarofsky, p. 109],

$$\frac{dT_e}{dt} = K(T_e)(T_i - T_e), \qquad (21)$$

where in general $K(T)$ is a slowly varying function of $T$. Equation (5.21) is usually represented in the following form, [Wood, p. 255],

$$\frac{dT_e}{dt} = \frac{2me}{m_i \tau_e}(T_i - T_e) \qquad (22)$$

where $\tau_e$ is the electron collision time interval, [Wood, p. 247].

In our reference frame equation (5.22) will be read as follow

$$\frac{\partial T_e}{\partial t} + u^i \frac{\partial T_e}{\partial x^i} = \frac{2m_e}{m_i \tau_e}(T_i - T_e),$$

which implies
$$u^1 \frac{dT_e}{dx} = \frac{2m_e}{m_i \tau_e}(T_i - T_e). \tag{23}$$

From [Wood, p. 247] we have $\tau_e^{-1} = f$, where $f$ is the collision frequency of each electron, thus from (5.23) we have
$$u^1 \frac{dT_e}{dx} = \frac{2m_e f}{m_i}(T_i - T_e). \tag{24}$$

Now it is sufficient to give another equation relating $T_e, T_i$ to $T$ and hence making it possible to determine $T_e, T_i$ when $T$ is known. Using (2.20), (2.51), (2.53) and (2.17) we have

$$T = \frac{1}{\nu_e + \nu_i}(\nu_e T_e + \nu_i T_i) = \frac{1}{(\rho_e/m_e) + (\rho_i/m_i)}\left((\frac{\rho_e}{m_e})T_e + (\frac{\rho_i}{m_i})T_i\right)$$
$$= \frac{1}{m_i \gamma_i - m_e \gamma_e}((m_i \gamma_i)T_e - (m_e \gamma_e)T_i), \tag{25}$$

If the plasma is singly ionized, i.e., if $e_i = -e_e = e$, the (2.25) simplifies to
$$T = \frac{1}{2}(T_e + T_i). \tag{26}$$

By the above argument we give equations for the mean variables as follows.

$$\rho u^1 = M,$$
$$-\tau^{11} + \rho(u^1)^2 - \frac{1}{2}\varepsilon(E^1)^2 + \frac{1}{2}\varepsilon(E^2)^2 - \frac{(B^1)^2}{2\mu} + \frac{(B^2)^2}{2\mu} + \frac{(B^3)^2}{2\mu} = \rho_1$$
$$-\tau^{12} + \rho u^1 u^2 - \varepsilon E_1 E_2 - \frac{1}{\mu}B^1 B^2 = P_2$$
$$-\tau^{13} + \rho u^1 u^3 - \varepsilon E_1 E_3 - \frac{1}{\mu}B^1 B^3 = P_3$$
$$\rho u^1(U_m + \phi + \frac{1}{2}(u)^2) + E_R - u^1 \tau^{11} - u^2 \tau^{21} - u^3 \tau^{31} = \frac{\rho u^1}{2\gamma_e \gamma_i \rho^2}J^2 - Q^1$$
$$= \frac{1}{\mu}(E^3 B^2 - E^2 B^3) + C \tag{27}$$

in which $M, P_1, P_2, P_3, C$ are constant. The required components of the stress tensor can be found from equation (2.34), (2.36). They are

$$\tau^{11} = -p_t + \tau_0^{11} = -p_t + \eta \frac{du^1}{dx}$$
$$\tau^{12} = \mu \frac{du^2}{dx}, \tau^{13} = \mu \frac{du^3}{dx} \tag{28}$$

in which we have set
$$\eta = \frac{4}{3}\mu + \zeta. \tag{29}$$

Only the $x$-component of the heat flux appears in the above equation. From equation (2.45) it is

$$Q^1 = K\frac{dT}{dx}, \tag{30}$$

To close our system of equations we need to use the generalized Ohm's law. If the momentum equation for the $\alpha$-species is multiplied by $\gamma_\alpha$, as done in the derivation of equation (5.13), the generalized Ohm's law is obtained. If the ions and electrons are not too far out of equilibrium with each other, the ion and electron stress tensors are of the same order of magnitude. Thus we have, up to $\dfrac{m_e}{m_i}$,

$$\gamma_e \tau_e^{ij} + \gamma_i \tau_i^{ij} = e\left(\frac{\tau_i^{ij}}{m_i} - \frac{\tau_e^{ij}}{m_e}\right) \simeq -e\frac{\tau_e^{ij}}{m_e} = \gamma_e \tau_e^{ij}, \tag{31}$$

we will further assume that the electron viscosity is negligible [An,Ak], that is we assume that

$$\tau_e^{ij} \simeq \delta^{ij} p_e. \tag{32}$$

Using (5.31) and (55.32) then (5.13) gives

$$\frac{1}{\gamma_e \gamma_i \rho}\frac{d}{dx}\left(-u^1\frac{dB^3}{dx}\right) = -E^2 - (u^1 B^3 - u^3 B^1)$$
$$+ \frac{\alpha_{12}}{\rho^2}\frac{(\gamma_e - \gamma_i)^2}{(\gamma_e \gamma_i)^2},\frac{dB^3}{dx},$$
$$\frac{1}{\gamma_e \gamma_i \rho}\frac{d}{dx}\left(-u^1\frac{dB^2}{dx}\right) = -E^3 - (u^1 B^2 - u^2 B^1) - \frac{\gamma_e + \gamma_i}{\gamma_e \gamma_i \rho}B^1\frac{dB^3}{dx}$$
$$- \frac{\alpha_{12}}{\rho^2}\frac{(\gamma_e - \gamma_i)^2}{(\gamma_e \gamma_i)^2},\frac{dB^2}{dx}. \tag{33}$$

Now using equations (5.27), (5.28), (5.30) and (5.33) we have the following equations,

$$\eta \frac{du^1}{dx} = \frac{1}{2\mu}((B^2)^2 + (B^3)^2) + Mu^1 - \left(P_1 + \frac{(B^1)^2}{2\mu}\right) + p_t,$$

$$\mu \frac{du^2}{dx} = Mu^2 \frac{1}{\mu} B^1 B^2 - (\varepsilon E_1 E_2 + P_2)$$

$$\mu \frac{du^3}{dx} = Mu^3 \frac{1}{\mu} B^1 B^3 - (\varepsilon E_1 E_3 + P_3) \tag{34}$$

$$\kappa \frac{dT}{dx} = M(U_m + \phi) + \frac{1}{2}M(u)^2 + E_R - u^1 \left(-p_t + \eta \frac{du^1}{dx}\right)$$

$$- u^2 \left(\mu \frac{du^2}{dx}\right) - u^3 \left(\mu \frac{du^3}{dx}\right) + \frac{1}{2\gamma_i \gamma_e M}(u^1 J)^2$$

$$- \frac{1}{\mu}(E^3 B^2 - E^2 B^3) - C,$$

$$\frac{u^1}{\gamma_i \gamma_e M} \frac{d}{dx}\left(-u^1 \frac{dB^3}{dx}\right) = -E^2 + (u^1 B^3 - u^3 B^1) + \frac{\gamma_e + \gamma_i}{\gamma_e \gamma_i M} u^1 B^1 \frac{dB^2}{dx}$$

$$+ \left(\frac{\alpha_{12}(\gamma_e - \gamma_i)^2}{\rho^2 (\gamma_e \gamma_i)^2}\right) \frac{dB^3}{dx},$$

$$\frac{u^1}{\gamma_i \gamma_e M} \frac{d}{dx}\left(-u^1 \frac{dB^2}{dx}\right) = -E^3 - (u^1 B^2 - u^2 B^1) + \frac{\gamma_e + \gamma_i}{\gamma_e \gamma_i M} u^1 B^1 \frac{dB^2}{dx}$$

$$- \left(\frac{\alpha_{12}(\gamma_e - \gamma_i)^2}{\rho^2 (\gamma_e \gamma_i)^2}\right) \frac{dB^3}{dx},$$

$$\frac{1}{\mu} \frac{dB_2}{dx} = J^3, \qquad \frac{dB^1}{dx} = 0,$$

$$\frac{1}{\mu} \frac{dB_3}{dx} = J^2, \qquad \frac{dE^1}{dx} = 0$$

In this equations we change the variables and constants as follow

$$u = u^1, v = u^2, w = u^3, P = \left(P_1 + \frac{(B^1)^2}{2\mu}\right),$$

$$P_2^* = (\varepsilon E^1 E^2 + P_2), P_3^* = (\varepsilon E^1 E^3 + P_3), \zeta_2 = u^1 J^3, \zeta_3 = -u^1 J^2$$

$$\zeta_2 = u^1 \frac{dB^2}{dx}, \zeta_3 = u^1 \frac{dB^3}{dx}, B_2 = B^2, B_3 = B^3, \tag{35}$$

$$\sigma = -\frac{(\gamma_e \gamma_i)^2}{(\gamma_e - \gamma_i)^2} \frac{\rho^2}{\alpha_{12}} = \frac{e^2 \nu_e}{m_e f_e}, \chi = \frac{\gamma_i + \gamma_e}{\gamma_e \gamma_i M} = \frac{m_i - m_e}{eM},$$

$$\beta = \frac{-1}{\gamma_e \gamma_i M} = \frac{m_1 m_2}{e^2 M}, E_1 = E^1, E_2 = E^2, E_3 = E^3$$

Then the set of equations (5.34) are written in the following form,

$$\eta \frac{du}{dx} = \frac{1}{2\mu}(B_2^2 + B_3^2) + Mu - P + p_t,$$

$$\mu \frac{dv}{dx} = mv - \frac{1}{\mu}B_1 B_2 - P_2^*,$$

$$\mu \frac{dw}{dx} = Mw - \frac{1}{\mu}B_1 B_3 - P_3^*,$$

$$\kappa \frac{dT}{dx} = M(U_m + \phi) + E_R - \frac{1}{2}M(u^2 + v^2 + w^2) - \frac{u}{2\mu}(B_2^2 + B_3^2)$$

$$- \frac{v}{\mu}B_1 B_2 - \frac{w}{\mu}B_1 B_3$$

$$Pu - P_2^* v - P_3^* w - \frac{1}{2}\beta(\zeta_2^2 + \zeta_3^2) - \frac{1}{\mu}(E, B_2 - E_2 B_3) - C, \qquad (36)$$

$$u \frac{dB_2}{dx} = \zeta_2,$$

$$u \frac{dB_3}{dx} = \zeta_3,$$

$$Bu \frac{d\zeta_3}{dx} = -E_2 + (uB_3 - wB_1) + \chi B_1 \zeta_2 - \sigma^{-1}\frac{\zeta_3}{u},$$

$$\beta u \frac{d\zeta_2}{dx} = E_3 + (uB_2 - vB_1) - \chi B_1 \zeta_3 - \sigma^{-1}\frac{\zeta_2}{u}.$$

In these equations $P_2^*$ and $P_3^*$ are constant and hence for the sake of simplicity we may set

$$P_2^* = P_3^* = 0,$$

and then we can write the equations in the following form,

$$\begin{aligned}
\mu \frac{dv}{dx} &= Mv - \frac{1}{\mu} B_1 B_2, \\
\mu \frac{dw}{dx} &= Mw - \frac{1}{\mu} B_1 B_3, \\
\eta \frac{du}{dx} &= \frac{1}{2\mu}(B_2^2 B_3^2) + Mu - P + p_t, \\
u \frac{dB_2}{dx} &= \zeta_2, \\
u \frac{dB_3}{dx} &= \zeta_3, \\
\beta u \frac{d\zeta}{dx} &= E_3 + (uB_2 - vB_1) - \chi B_1 \zeta_3 - \sigma^{-1} \frac{\zeta_2}{u}, \\
\beta u \frac{d\zeta}{dx} &= -E_2 + (uB_3 - wB_1) + \chi B_1 \zeta_2 - \sigma^{-1} \frac{\zeta_3}{u}, \\
\kappa \frac{dT}{dx} &= M\left(U_m + \phi + \frac{E_R}{M}\right) - \frac{1}{2}M(u^2 + v^2 + w^2) \\
&\quad - \frac{u}{2\mu}(B_2^2 + B_3^2) - \frac{v}{\mu} B_1 B_2 - \frac{w}{\mu} B_1 B_3 \\
&\quad - Pu - \frac{1}{2}\beta(\zeta_2^2 + \zeta_3^2) - \frac{1}{\mu}(E_3 E_2 - E_2 B_3) - C.
\end{aligned} \qquad (37)$$

These are the equations for the shock layer equations for a fully ionized and singly ionized plasma. Using these equation and the equation (3.15)-(3.30) one can determine the partial variables when the mean variables are determined. As we have seen during the derivation of (5.37), the variables and constants in (5.37) are as follow, mean velocity vector $= (u, v, w)$,

Magnetic field $= (B_1, B_2, B_3)$,

Electric field $= (E_1, E_2, E_3)$,

mean density $= \rho = Mu^{-1}$,

Specific volume $= \rho^{-1} = \frac{u}{M}$,

electric conductivity $= \sigma$,

collision frequency of electrons $= f_1$,

Hall effect coefficient $= \chi$,

second viscosity coefficient $= \mu$,

first viscosity coefficient $= \eta$,

heat conduction coefficient $= \kappa$,

heat flux vector $= (Q^1, Q^1, Q^3)$,

current density $= (J^1, J^2, J^3)$.